\documentstyle[12pt]{article}

\newcounter{const}

\begin{document}
%\large
%\normalsize
% 26.11.2004
\begin{center}

{\bf CORRECT SOLVABILITY OF NONLINEAR ORDINARY DIFFERENTIAL EQUATIONS
IN ORLICZ SPACES.}\\
\vspace{1mm}
{\bf E.I.OSTROVSKY.} \footnote{ Department of Mathematics, 
Ben - Gurion University, Beer - Sheva, Israel.\\
e - mail: galaostr@cs.bgu.ac.il} \\
\vspace{1mm}

\end{center}

{\bf Abstract.} {\it We prove in this article the well posedness of  non - linear Ordinary 
Differential Equations (ODE) of first and second order in Orlicz spaces  with unbounded 
domain of definition.}\\ 

\vspace{1mm}

{\bf Key Words.} Orlicz spaces, ordinary non - linear differential equations, equivalent norms.\\

\vspace{1mm}

{\bf AMS (MOS) subject classification. } Primary 34C11, 34B40, 47E05.\\

\vspace{2mm}

{\bf 0. Introduction. Statement of problem.}\\

 Let us consider the first order  non - linear ODE of a view:
$$
  dy/dx - q(y) = g(x), \ \ x \in R. \leqno(1)
$$
or (non - linear Sturm - Liouville equation)

$$
 d^2y/dx^2 - q(x,y) = g(x), \ x \in R. \leqno(2)
$$

 It is proved in the works ([2]; [3], [4], [5]) that under some conditions (necessary 
conditions and sufficient conditions) the  equations (1) and  (2) {\it without boundary conditions}
 are correct solvable in the spaces 
$ L_p(R) $ and consequently, under some simple additional conditions
in the correspondent Sobolev spaces $ W_1^p(R), W_2^p(R). $ \par
 {\bf Our goal is some generalizations of those results  on the Orlicz spaces $ L(N) $
with  $ N \ - $ Orlicz function $ N = N(u) $ instead classical functions } $ |u|^p, \ 
p = const \ge 1, $  in particular, on the Orlicz  spaces $ L(N) $ 
{\it without } the so - called $ \Delta_2 $ condition, e.g. Exponential Orlicz Spaces 
(EOS). This allow us to find some {\it new} properties of solutions, for instance, to prove 
the exponential integrability of solutions and its derivatives. \par
 Probably, it is very interest to describe {\it all } the Orlicz spaces for which the equations
 (1), (2) are correct solvable.  The statement of this problem belongs to L.Shuster. 
 Now this general problem is open, but we can prove the correct solvability of our 
equations on two important classes of Orlicz spaces.\par
 Recall here that the equation (non - linear, in general case)
$ Ay = g $ is called correct solvable (or, in other hand, well posed)
(more exactly, Lipshitz correct solvable) on the Banach space 
$ B_1 $ into the space $ B_2, $ if for all $ g \in B_1 $ the solution $ y = A^{-1}g $ 
there exists, is unique, belongs to the space $ B_2 $ and 

$$
||A^{-1} g||B_2 \le C ||g||B_1; \ ||A^{-1} g_1 - A^{-1}g_2||B_2 \le C ||g_1 - g_2||B_1,
$$
where $ C = const $ does not depend on $ g, g_1, g_2. $ \par

\vspace{2mm}

{\bf 1. Description of using Orlicz spaces.}\par
 We will consider a two kinds of Orlicz spaces on the real line $ R $ with usually 
(unbounded) Lebesque measure $ \mu. $  Recall here  that  if the function 
$ N = N(u),  u \in R^1  $ is some $ N \ - $ Orlicz function (even, 
downwards convex, $ N(u) \ge 0, \ N(u) = 0 \ \Leftrightarrow u = 0, $
strong increasing in the self - line  $ R^1_+, $ etc.), then the Orlicz norm $ ||f||L(N) $
of a (measurable) function $ f: R \to R $ relative to the $ N \ - $ Orlicz function 
$ N = N(u) $ may be defined by the formula
$$
||f||L(N) = \inf \{k, \ k > 0, \ I(N(|f|/k)) \le 1 \}.
$$
 In this paper  $ I(f) = \int_R f(x) \ dx. $ As a particular case, if $ N(u) = 
|u|^p, \ p = const \ge 1 $
we obtain the classical $ L_p = L_p(R) \ $ spaces with the norm
$$
|f|_p \stackrel{def}{=} I^{1/p} |f|^p.
$$ 

 {\bf A}. We define a class $ \Delta $ as a {\it set} of all Orlicz spaces with correspondent 
$ N \ - $ function $ N = N(u) $ belonging to all the classes $ \Delta_2(\infty), 
\nabla_2(\infty)  $ in the terminology of the book [1], p. 22 - 24.
 By definition,
$$
\Delta_2(\infty) = \{N:\exists k < \infty, \forall u \ge 0 \ \Rightarrow  N(2u)/N(u) <k \};
$$
$$
\nabla_2(\infty) = \{N: \ \ \exists l > 1, \forall u > 0 \ \Rightarrow N(2u) \le  N(l u)/(2l) \}.
$$

or, briefly, $ N(\cdot) \in \Delta \stackrel{def}{=} \Delta_2(\infty) \cap \nabla_2(\infty). $ 
For example, let $ N(u) = N(m,r; u) \stackrel{def}{=} 
 |u|^m \ \left(\log^r( \exp(m + |r|) + |u|) \right), 
m = const > 1, r = const \in R, $ then $ N(\cdot) $ is some $ N \ - $ Orlicz function such that
$ N(\cdot) \in \Delta. $ \par
{\bf B.} An other very important class of $ N \ - $ Orlicz functions are so - called 
$ EOF \ = $ Exponential Orlicz Functions  and correspondent Exponential Orlicz Spaces $ EOS $ 
will be considered. Let 
$ \varphi = \varphi(z), \ z \ge 1 $ be some continuous  function such that 
the function $ h(y) = h_{\varphi}(y):= \varphi(\exp y), \ y \ge 0 $ is strong increasing, downward 
convex and 
$$
\sum_{k=3}^{\infty} \exp(h(k) - h(k+1)) < \infty. \leqno(3)
$$
 The set of all those function we will denote $ \Phi; \ \Phi = \{\varphi \}. $

  For example, put $ \varphi(z) = \varphi_{m,r}(z) = z^m \  \log^r [(\exp(m + |r|) + z], $
 but here $ m = const \in (0, \infty); $  
or $ \varphi(z) = \varphi_{\beta}(z) = \log^{1+ \beta}( 2 + z), \ \beta = const > 0.  $ 
Then $  \varphi_{m,r}(\cdot) \in \Phi, \ \varphi_{\beta} (\cdot) \in \Phi. $ \par

 Let also $ \alpha = const \ge 1, \ \varphi \in \Phi. $ We denote by 
$ \exp_{\alpha}\varphi(z) $  the following 
continuous  function: at  $ z \in [0,C_1] \ \Rightarrow  
\exp_{\alpha} \varphi(z) = C_2 z^{\alpha}, $ and at 
$ \ z\in (C_1, \ \infty) \ \Rightarrow  \exp_{\alpha} \varphi(z) = \exp \varphi(z).$ \par
 Let us prove at first the existence of the constants 
 $ C_1 = C_1(\alpha, \phi), C_2 = C_2(\alpha, \phi) $ {\it 
such that} $ \exp_{\alpha}(\cdot) \in EOF. $ We will use 
 the so - called Young - Fenchel, or Legendre transform:

$$
h^*(w) = \sup_{y \ge C} (yw - h(y)) = y_0 w - h(y_0),  \ \ w \ge C_1 = C_1(C),
$$
 the value $ y_0 = y_0(w) $ there exists and is unique for all sufficiently larges values $ w.$ \par
 We must only prove  that the constant $ C_2 $ there exists and is nontrivial:
$$
C_2 = \inf_{z > C_1} \exp \varphi(z)/z^{\alpha} = \exp \left(-\sup_{z \ge C_1} (\alpha \log z -
\varphi(z) \right) =
$$

$$
 \exp \left( -\sup_{y \ge \exp(C_1) } (\alpha y - h(y) \right) = \exp \left(- h^*(\alpha) \right);
$$
therefore the constant  $  C_2 = C_2(\alpha, \varphi(\cdot) ) $ there exists and is nontrivial:
$ C_2 \  \in (0, \infty). $ \par
 {\it Further we will choose the constants $ C_1, C_2 $ only such that the function }
$ \exp_{\alpha} \varphi(\cdot) \in EOF. $ \par
{\bf Definition.}  We define the so - called Exponential Orlicz Space (EOS) 
$ B(\alpha; \varphi), $ where $ \ \alpha = const \ge 1, \ \varphi(\cdot) \in \Phi $
 as the Orlicz space on the set $ R^1 $
with the correspondent $ N \ - $ Orlicz function
$$
N(u) = N(\alpha; \varphi, u) = \exp_{\alpha} \varphi(u)
$$

and the correspondent Orlicz norm $ ||f||L(\exp_{\alpha}(\varphi)) = ||f||B(\alpha; \varphi) $
and will say this spaces as EOS = Exponential Orlicz Spaces. \par
 In the case $ \varphi(z) = \varphi_m(z) = z^m, \ z,m >0 $ we will write simply 
$ ||f||B(\alpha, m) = ||f||B(\alpha; \varphi_m), $ i.e. here 
$$
 N(z) = N(\alpha,m, z) = \exp_{\alpha}(z^m), \ z \ge 0. 
$$ 

 Note that at $ \alpha = m(l+1), \ l = 0,1,2,\ldots $ the equivalent $ N \ - $ function for
$ N(\alpha,m;u) $ may be constructed by the formula ( see [9], [10], p. 13 - 14)
$$
N(\alpha,m;u) \sim \exp \left(|u|^m \right) - \sum_{k=0}^l |u|^{mk}/k!.
$$
 In the case $ m = \infty $ the space $ B(\alpha, \infty) $ may be defined as a projective
limit of the  spaces $ B(\alpha,m) $ at $ m \to \infty. $ But $ B(\alpha, \infty) $ 
 is isomorphic to the space of all 
bounded $ (mod \ \ \mu) \ $ and integrable with power  $ \alpha $ functions: 
$ ||f||B(\alpha,\infty) \le $
$$
C_1(\alpha) \left[ \ |f|_{\alpha} + vraimax_{x \in R} |f(x)| \right] 
\le C_2(\alpha) ||f||B(\alpha, \infty).
$$

\vspace{2 mm}

{\bf 2. Main results.} We consider at first  the  differential equation (1),
or, by notation, $ y = Q[g], $ if obviously the solution $ y $ there exists and is unique.
 We suppose  $ q(0) = 0 $ and 
$$
0 < m = m(q) \stackrel{def}{=} \inf_{x \ne z} |q(x) - q(z)|/|x-z| \le 
$$
$$
 \sup_{x \ne z}|q(x) - q(z)|/|x-z| \stackrel{def}{=} M(q) = M < \infty.\leqno(4)
$$ 
Recall also here the definition of Qrlicz - Sobolev norms and correspondent Orlicz - Sobolev
 spaces: $ W_k(L(N)), \ k = 1,2,\ldots $ consists on all the measurable functions $ y = y(x) $
with finite norm 
$$
||y||W_k(L(N)) = \sum_{l=1}^k ||d^l y/dx^l||L(N) + ||y||L(N).
$$

{\bf Theorem 1.} {\it Suppose that } $ N(\cdot) \in \Delta $ {\it or } $ N(\cdot) =  
B(\alpha; \varphi) $ {\it for some } $ \alpha \ge 1, \ \varphi \in \Phi. $
{\it Then the problem (1) is Lipshitz correct solvable on the space $ L(N) $ into the space }
$ W_1(L(N)): $

$$
 ||Q[g]||W_1(L[N]) \le C||g||L[N], \ C = C(m,M,N(\cdot)) \in  (0,\infty);
$$
$$
 ||Q[g_1]- Q[q_2]||W_1(L[N]) \le C(m,M, N(\cdot)) \ ||g_1-g_2||L[N]. \leqno(5)
$$
 Note that the assertion (5) is some generalization of main result of paper [2]. It is 
proved at the same place that the conditions  (4) are necessary and sufficient for (5) even 
for the spaces $ L_p. $  \par 
  We consider now the problem 2.  Denote the solution of equation (2) 
by $ y = S[g] $ again in the Orlicz space $ L[N](R) = \{g\}. $  (We will prove further the
existence and uniques of $  S[g].) $ \par
 Let us introduce for the finite measurable function $ v = v(x,y), \ x,y \in R $ of two variables
with condition $ v(x,0) = 0 $  the following norm:

$$
|||v||| = \sup_x \ \sup_{y \ne z} |v(x,y) - v(x,z)|/|y-z|.
$$
 We denote by $ V $ the Banach space of all the (measurable) functions $ v $ such that $ v(x,0) = 0 $
with finite norm $ |||\cdot|||: \ \ V = \{v: \  |||v||| < \infty. \} $
 We suppose that the measurable function $ q(x,y) $ has a properties: $ 1) \ q(x,y) \ge 0; $ 
$ 2) \ \ \exists v \in V, \exists q_0(\cdot) \in L_1^{loc},  \ \ q_0 \ge 0, \ v(x,0) = 0, $
$$
q(x,y) = q_0(x)y + v(x,y).
$$
 
3) We  define also 
 $ d = d(x) $ as a (unique)  non - negative  solution of equation
$$
 \int_0^{d \sqrt{2}} \ dt \ \int_{x-t}^{x+t}  q_0(\xi) d \xi = 2, \leqno(6)
$$
 and denote 
$$
A = \inf_{x \in R} d(x), \ \ B = \sup_{x \in R} d(x), \leqno(7)
$$

$$
\nu(B) = 8 \exp(-1/e)  B \ \max(B,1).
$$

{\bf Theorem 2.} {\it Suppose $ B < \infty $ and

$$ 
|||v||| < 1/\nu(B). \leqno(8)
$$
  Assume again  that  $ N(\cdot) \in \Delta  $   or  $ N(\cdot) = B(\alpha; \varphi) $
 for  some  $ \alpha \ge 1, \ \varphi \in \Phi. $  We assert that  for all functions $ N(\cdot) $  
and potentials  $ \ q(x,y) $  which satisfies our conditions 
the problem (2) is well - posed on the space $ L_{\beta} = L_{\beta}(R),$ where $ \ \beta \ge 
\ \alpha, $ into the space $ L[N]: $
$$
||S[g] ||L[N] \le C_4(q(\cdot)) \ ||g||L_{\beta}, 
$$
$$
||S[g_1] - S[g_2]||L[N] \le C_4 \ ||g_1-g_2||L_{\beta}. \leqno(9)
$$
  Suppose in addition  that

$$
\sup_x q_0(x) = C_5 < \infty.
$$

 Then the problem 2 is correct solvable in the space $ L_{\beta} $  into the space }
$ W_2(L(N)): $

$$
||S[g]||W_2(L(N)) \le C_4(q(\cdot)) \ ||g||L_{\beta}, 
$$

$$
||S[g_1] - S[g_2]||W_2(L(N)) \le C_4 \ ||g_1-g_2||L_{\beta}.
$$

{\bf Theorem 3.} {\it Suppose that in the problem (2) $  v(x,y) = 0 $ (linear equation) 
and }

$$
A \stackrel{def}{=} \inf_x d(x) > 0. \leqno(10)
$$
{\it  Then  for all $ \alpha > 1, \ \beta > \alpha, \ \delta \in (0,\beta -\alpha) $  the
 problem (2) is ill - posed in the space $ L_{\beta - \delta} $ into the space } 
$ L_{\beta}. $ Namely, $  \ \forall \beta > \alpha, \ 
\delta \in (0, \beta -\alpha) \ \  \exists g(\cdot) \in L_{\beta- \delta} \ \Rightarrow  S[g] 
\notin L_{\beta}. $ \par
  Theorem 2 is some generalization of main result of paper [3]. It is obtained in [4]  
in linear case $ v = 0 $ the criterion of correct solvability (2) in the spaces $ L_p 
\to L_p, \ p \ge 1. $  \par

 {\bf Remark 1.} We can notice the diffrence between equations of first order ODE (1)  and second 
order  ODE (2).  In first case the right - side of equation must belong to the Orlicz space 
$ L(N), $ in the second case $ g(\cdot) $ must belong only to the $ L_{\beta}(R) $ space.\par

\vspace{2mm}

{\bf 3. Auxiliary result.}  Denote $ \psi(p) = \psi(p; \varphi) = 
\exp \left(h^*(p)/p  \right). $  Let us introduce 
a {\it new } Banach space $ G(\alpha; \varphi), \ \alpha \ge 1, \ \varphi \in \Phi, $ 
as a set of all measurable functions $ f:R \to R  $ with finite norm

$$
||f||G(\alpha, \varphi) \stackrel{def}{=} \sup_{p \ge \alpha} |f|_p /\psi(p) < \infty.
$$

{\bf Theorem 4.} 
 {\it We propose that the norms } $ || \cdot ||B(\alpha, \varphi) $ {\it and } 
 $ || \cdot ||G(\alpha; \varphi) $
{\it are equivalent: } $ \exists C_3, C_4 = C_3,C_4(\alpha, \varphi) \in (0,\infty) \ 
\Rightarrow $

$$
C_3 \ ||f||G(\alpha; \varphi) \le ||f||B(\alpha; \varphi) \le C_4 \ ||f||G(\alpha; \varphi).\leqno(11)
$$

{\bf Proof of theorem 4.} 
 Assume at first that $  ||f||B(\alpha; \varphi) < \infty. $
Without loss of generality we can suppose  

$$
I \left(\exp_{\alpha} \varphi(|f|) \right) = 1. 
$$
 Let us introduce the function

$$
\gamma(p) = \gamma_{\alpha}(p) = \sup_{z > 0} z^p / \exp_{\alpha} \varphi(z).
$$
 We have for the values $ p \ge \alpha $ and some $ C_1 = C_1(\alpha, \varphi(\cdot)) \in 
(0, \infty): $ 

$$
\gamma(p) \le  \max \left[\max_{z \in (0, C_1] } C^p  z^{p-\alpha},  \ \sup_{z \ge C_1} z^p \ 
\exp (- \varphi(z)) \right]  \le
$$

$$
\max \left[ C_2^{p-\alpha}, \ \exp( \sup_{z \ge C_1} (p \log z - \varphi(z))) \right] =
$$
$$
 \max \left[ C^{p-\alpha}, \ \exp(\sup_{v \ge C_3} (p v - h(y))) \right] =
$$

$$
\max \left[C^{p-\alpha}, \  \exp h^*(p) \right] \le C_4^p(\alpha) \exp h^*(p).
$$

 Following, for the values $ p \ge \alpha $ we have: $ z \ge 0 \ \Rightarrow $

$$
z^p \le \gamma(p) \ \exp_{\alpha} \varphi(z) \le C_4^p(\alpha) \ \psi^p(p) \ \exp_{\alpha}
\varphi(z). 
$$

Therefore
$$
|f|^p \le C^p(\alpha) \ \psi^p(p) \ \exp_{\alpha}( h^*(|f|)), \ \ 
|f|_p \le C(\alpha; \varphi) \ \psi(p),
$$

$$
||f||G(\alpha, \varphi) \le C(\alpha; \varphi(\cdot)) < \infty.
$$

 Inverse, assume that 

$$
|f|_p^p  \le \exp \left(h^*(p) \right), \ p \ge \alpha.
$$
 We have by virtue of Chebyshev inequality for all the values $ w \ge C_5:$

$$
T(|f|,w) \stackrel{def}{=} \mu \ \{x: |f(x)| > w \} \le 
\exp \left(h^*(p) - p \log w \right),
$$
 After the minimization of the  right - side over $ p, \ p \ge C, $ we receive 
for $ w \ge C_2: $

$$
T(|f|,w) \le \exp \left(- h^{**}(\log w) \right) = \exp \left(- h(\log w) \right)
$$
  on the basis of theorem of Fenchel - Moraux. We conclude for  the 
 value of $ \varepsilon = \exp(-2), $ choosing $ W(k) = \exp(k) $ and denoting

$$
U(k) = U(|f|,k) = \{x: \ W(k) \le |f(x)| < W(k+1) \}: 
$$

$$
I(\exp(\varphi_{\alpha}(\varepsilon |f|))) \le C + 
\sum_{k=3}^{\infty} \int_{U(k)}\exp(\varphi(\varepsilon |f|) ) \ dx \le 
$$

$$
C+\sum_{k=3}^{\infty} \exp \left[\left( h(\varepsilon \  W(k+1) \right) \ \cdot \ 
T(|f|, W(k)) \right] \le
$$

$$
C + \sum_{k=3}^{\infty} \exp(h(k) - h(k+1)) < \infty
$$
by virtue of condition 3. This completes the proof of theorem 4.\par

For example, let $ N(u) = N_{\alpha,m}(u)= \exp_{\alpha}(\varphi_m(u)) = 
\exp_{\alpha} |u|^m. $  It follows from theorem 4 that 
$$
||f||L(N_{\alpha,m})  < \infty \ \Longleftrightarrow \ \sup_{p \ge \alpha} |f|_p \ 
p^{-1/m} < \infty,
$$
 or equally

$$
\exists \varepsilon > 0, \ I(\exp_{\alpha}(\varepsilon |f|) < \infty \ \Longleftrightarrow
\ \sup_{p \ge \alpha} |f|_p \ p^{-1/m} < \infty.
$$
{\bf Notice.} Let us introduce the {\it weight Lorentz } norm:

$$
||f||^*_bG(\alpha,\varphi) = \sup_{p \ge \alpha} ||f||_{p,b}/\psi(p),
$$
where $ ||f||_{p,b} $ is the Lorentz norm (more exactly, seminorm):

$$
||f||_{p,b} = \left[\int_0^{\infty} T^{p/b}(|f|,x) \ dx^b \right]^{1/b}, 
$$
$ p \in [1, \infty), \ b \in [1,\infty], $ where if $ b = \infty $ then
$$ 
||f||_{p,\infty} = \sup_{x \ge 0} \left( x \ T^{1/p}(|f|,x) \right).
$$
 
 It is easy to prove using the embedding theorem for the Lorentz spaces  as well as by 
proving of theorem 4 that all the norms 
$$
||\cdot||B(\alpha; \varphi), \ ||\cdot|| G(\alpha; \varphi), \ ||\cdot||^*_b G(\alpha, 
\varphi)
$$ 
are equivalent with constants does not depending on $ b. $ \par
 Note than if we consider the Orlicz space $ L(N), \ N \in EOF $ on the arbitrary 
measurable space $ (\Omega, F, \mu) $ with finite measure $ \mu, $ the result of theorem 4 
is known (see [12], p.341). \par

\vspace{2mm}
{\bf Proof of theorem 1. } Let us consider at first the case $ N \in \Delta.$ Let 
$ g \in L(N), \ N \in \Delta.$  We will use the main result of paper [2]:

$$
\exists \ C(m,M) \in (0,\infty), \ \forall p \ge 1 \ |Q[g]|_p \le C(m,M) \ |g|_p, \leqno(12)
$$ 
and 
$$
|Q[g_1] - Q[g_2]|_p \le C(m,M) \ |g_1 - g_2|_p,
$$
where we denote for this problem $ y = Q[g]. $ It follows from (12) that the operator $ Q $ is 
correct defined  and bounded as operator $ L_p \to L_p, \ p \ge 1. $ The first proposition 
of theorem 1 follows from Ryan's theorem ( [1], p. 193).\par 
 Let now $ N(\cdot) = B(\alpha; \varphi) $ for some $ \alpha \ge 1,\ \varphi \in \Phi $
and let $ g(\cdot) \in B(\alpha; \varphi). $
By virtue of theorem 4 $ |g|_p \le C_1 \ \psi(p)||g||B(\alpha; \varphi), \ p \ge \alpha. $ 
Therefore (see (12) )

$$
|Q[g]|_p = |y|_p \le C_1 \ C(m,M) \ ||g||B(\alpha; \ \varphi) \ \psi(p).
$$
 Again from theorem 4 follows 
$$
||Q[g]||B(\alpha, \varphi) = ||y||B(\alpha;\varphi) \le C_1 \ C_2 \ C(m,M) \ 
||g||B(\alpha; \varphi),
$$
and analogously 
$$
||Q[g_1] - Q[g_2]||B(\alpha; \varphi) \le C_3(\alpha, \varphi, m,M) 
 \ ||g_1 - g_2||B(\alpha; \varphi).
$$
 This completes the proof of theorem 1.\par
{\bf Proof of theorem 2.} {\it Part 1.} We consider here the linear case, i.e. 
$ v(x,y) = 0.$ We can denote by $  y_0 = S[g] $ the solution of linear equation
$$
d^2y_0(x)/dx^2 - q_0(x) \ y_0(x) = g(x), \ \ \lim_{|x| \to \infty} y_0(x) = 0,
$$
as long as $ y_0 $ there exists and is unique ([4], [5]).
 The first part of this theorem is proved analogously to 
the proof of theorem 1, since (see [4]) 

$$
|S[g]|_p \le C \ |g|_p, \ \ |S[g_1] - S[g_2]|_p \le C \ |g_1 - g_2|_p,
$$
{\it Part 2.}  Further, we will denote by $  \Gamma(t,x) $ the Green's function for 
the linear equation (2):

$$
y_0(x) = \int_R \Gamma(x,t) \ g(t) \ dt.
$$
 We will use the fine result of paper [5]: the function $ \Gamma(\cdot, \cdot) $ there exists, 
is unique, and
$$
\Gamma(x,t) = \sqrt{\rho(x) \rho(t)} \ \exp \left( - 0.5 \left| \int_x^t d\xi/\rho(\xi) \right|
\right),
$$
where 
$$
2^{-3/2} d(x) \le \rho(x) \le 2^{-1/2} d(x). \leqno(13)
$$
 It follows from (13) and condition (8) 

$$
\Gamma(x,t) \le B \ \sqrt{2} \ \exp \left( -2^{-3/2} |t-x|/B \right). \leqno(14)
$$
 Therefore
$$
|y_0(x)| \le B \ \sqrt{2} \ \int_R \exp \left( - 2^{-3/2} |t-x|/B \right) \ g(t) \ dt =
$$
$$
s*g(x), \ \ \ s(x) = B \ \sqrt{2} \ \exp \left(-2^{-3/2} |x|/B \right),
$$
and the symbol $ f*g $ denotes the usually convolution for the function defined on $ R. $
Using the Young inequality for convolution we obtain for all the values $ r \ge \beta:$
$$
|S[g]|_r \le |g|_{\beta} \cdot \sup_{p \in [1, \beta/(\beta-1)) } \ |s|_p,
$$
where at $ \beta = 1 \ \Rightarrow \beta/(\beta-1) = + \infty.$ \par

 It is easy to calculate that 
$$
\sup_{p \ge 1} \ |s|_p \le 8 \ \exp(-1/e) B  \max(B,1) = \nu(B) < \infty,
$$
following, 
$$
\sup_{p \ge \alpha} | \ S[g] \ |_p \le \nu(B) \ |g|_{\beta}, 
$$
and $ \forall \varphi \in \Phi, \ \alpha \ge 1, \ \beta \ge \alpha $

$$
||S[g]||G(\alpha; \varphi) \le C(\alpha; \varphi) \ \nu(B) \ |g|_{\beta}. \leqno(15)
$$
{\it Part 3.} Let us consider in this section the non - linear case $ v \ne 0. $ We can 
rewrite the equation (2) on the form

$$
y = W[y],  \ \ W[y](x) = S[g](x) + \int_R \Gamma(x,z) \ v(z,y(z)) \ dz. \leqno(16)
$$
Since $ |v(x,y)| \le C |y|, $ it is evident that the (non - linear) operator $ W[\cdot] $
has the property: $ W: L_{\beta} \to L_{\beta}, $ i.e. 
$$
||W[\cdot]||(L_{\beta} \to L_{\beta} ) < \infty
$$
and we have by virtue of inequality (14):

$$
W[y_1](x) - W[y_2](x) = \int_R \Gamma(x,z) [v(z,y_1(z) - v(z,y_2(z)] \ dz;
$$

$$
|W[y_1](x) - W[y_2](x)| \le |||v||| \ \int_R s(x-z) |y_1(z) - y_2(z)| \ dz;
$$

$$
|W[y_1] - W[y_2]|_{\beta} \le |||v||| \ \nu(B) \ |y_1 - y_2|_{\beta}.
$$
 Therefore, the operator $ W[\cdot] $ satisfies the contraction property in the space 
$ L_{\beta}. $  Following, there exists the fixed point of $ W[\cdot] $ in the space 
$ L_{\beta}. $ \par
 The statement of part 3  follows from  part 1 and 
(16), as long as all the terms of
the right side (16)  belong to the space $ B(\alpha; \varphi); \ \alpha \ge 1, 
\varphi \in \Phi.$ \par
{\it Part 4.} Suppose in addition  $ \sup_{x \in R} q_0(x) = C_6 < \infty, \ 
v(x,y) = 0. $ We can rewrite in this case our equation (2) on the form

$$
d^2y/dx^2 = g(x) + q_0(x) y(x),
$$
in the spaces $ L_{\beta} \to L_{\beta}, $ i.e. we assume that 
$ g \in L_{\beta}, \ $ and following, $
y \in L_{\beta}. $ Therefore  $ d^2y/dx^2 \in L_{\beta}. $
The last assertion of theorem 2  follows from the classical Kolmogorov's inequality

$$
||dy/dx||_p^2 \le 16 ||y||_p \ ||d^2y/dx^2||_p;
$$
 see, for example, [11], p. 49. This completes the proof of theorem 3. \par

{\bf 4. Proof  of  theorem 3.} Let in the equation (2) $ \alpha > 1, v(x,y) = 0. $  
It follows from (13) that
$$
\Gamma(x,t) \ge 2^{-3/2} A \ \exp \left(-2^{-5/2} |t-x|/A \right). \leqno(16)
$$
 Let $ g(x) = g_{\beta}(x) = (x \log^2 x)^{1/\beta} $ if $ x \ge 2 $ and $ g(x) = 0, \ 
x < 2. $ Then $ g \in L_{\beta} $ and $ g \notin L_{\beta - \delta} \ \forall \delta 
\in (0, \beta - \alpha). $  We have for the values $ x \ge 3: $
$$
W[g] \ge C_7 \int_2^{\infty} \exp(-|x-t|) \ (t \ \log^2 t)^{1/\beta} \ dt = 
$$

$$
= C_7 x^{1-1/\beta} \ \int_{2/x}^{\infty} \exp(-x|z-1|) \ (z \ \log^2(xz))^{-1/\beta} \ dz.
$$
 The exact asymptotic of  the last integral at $ x \to \infty  $
may be calculated by means of the Laplace's method; the critical point $ z = 1. $ We obtain:

$$
W[g](x) \sim C_8 (x \ \log^2 x)^{-1/\beta}, \ \ C_8 = C_8(\beta) \in (0,\infty).
$$
 Therefore, $ W[g] \in L_{\beta} $ but $ \forall \delta \in (0, \beta - \alpha) \ \Rightarrow 
W[g] \notin L_{\beta - \delta}. $ \par

Finally, let us consider in addition to the problem (1) and (2) the 
{\it Pseudodifferential} \ linear operator $ P $ in the space $ R^n $ with the symbol $ P(\xi). $
Assume that $ \forall k = 1,2,\ldots, [n/2]+1, \ ([x] $ denotes here the integer part of $ x) $  

$$
\sup_{z > 0} \max_{|\zeta | \le k} z^{-n} \int_{\xi: z < |\xi| \le 2 z}  |\xi|^{k} \ 
|P^{(\zeta)}(\xi)|^2 \ d \xi < \infty;
$$
 where  $ \xi = (\xi_1,\xi_2,\ldots,\xi_n), \
 |\xi| = \sum_i |\xi_i|, \ \zeta = (\zeta_1,\zeta_2, \ldots, \zeta_n), \
\zeta_j = 0,1,\ldots;  $
$$
P^{(\zeta)}(\xi) = \frac{\partial^{(|\zeta|) } P(\xi)}{\prod_{i=1}^n \partial^{\zeta_i} \xi_i}.
$$
 It is known (see, for example, [13], p. 262 - 270)) that 
$$
|Pf|L_p(R^n) \le C(n,P) \ p \ |f|L_p(R^n), \ p > 1. \leqno(17)
$$
 Suppose  $ f \in B(\alpha,m) $ for some $ \alpha > 1, \ m > 0. $ From theorem 4 follows that
$ Pf \in B(\alpha, m/(m+1)), $ hence  $ P: B(\alpha,m) \to B(\alpha, m/(m+1)) $ and 
$$
|| P ||( B(\alpha,m) \to B(\alpha, m/(m+1))) < \infty.
$$

\vspace{2mm}
{\bf Aknowledgements}. I am very gratefull to prof. N.Chernyavskaya and prof. L. Shuster 
for many fruitful  consultations about considering here problems. \\ 
\vspace{2mm}

{\bf Concluding remark.} Our results (without proof) was announced in [14].\\

\vspace{3mm}
{\bf References.} \\

1. M.M.Rao, Z.D.Ren.  Applications of Orlicz Spaces.{\it  Marcel Dekker Inc.,} 2002. New York,
Basel.\\
2. V.E. Slyusarchuk.  Necessary and Sufficient Conditions for the Lipshitzian Invertibility 
of the Nonlinear Differential Mapping $ d/dt - f $ in the Spaces $ L_p(R,R), 1 \le p \le \infty. $
{\it Math. Notes,} 2003, {\bf 73}, $ N^o $ 6, 843 - 854.\\
3. N.Chernyavskaya, L.Shuster.  Weight Summability of Solution of the Sturm - Liouville 
Equation. {\it Journal of Diff. Equations,} {\bf 151,} 1999,  456 - 473. \\
4. N.Chernyavskaya and L.Shuster. A Criterion for Correct Solvability of the 
Schturm - Liouville Equation in the Space $ L_p(R). $ {\it Proceedings of the American Mathematical 
Society,} 2001, V. 130 Number 4, p. 1043 - 1054.\\
5. N.Chernyavskaya and L.Shuster. Estimates for the Green Function of a General Sturm - 
Liouville Operator. {\it J. of Differential Equations}, 1994, {\bf 111}, p. 410 - 420.\\
6. E.B.Davies and E.M.Harell.  Conformally flat Riemannian metrics, Schrodinger 
operator and semiclassical approximation. {\it J. Diff. Eq.,} 1987, {\bf 66}, p. 165 - 188. \\
7. D.Medvedev and V.Vlasov. On certain properties of exponential solutions of difference 
differential equations in Sobolev spaces. {\it Functional Differential Equations.}
 V.9, $ N^o $ 3 - 4, (2002), 423 - 435.\\
8. Abdelnaser J. Al - Hasan. \ A Note on a maximal singular integral. {\it Functional Differential 
Equations.} V.5, $ N^o $ 3 - 4, (1998), 309 - 314.\\ 
9. N. Trudinger. On imbeddings into Orlicz Spaces and some applications. {\it J. Math. Mech.,}
1967, v. 17,  \ 473 - 483. \\
10. M.E.Taylor. Partial Differential Equations, III (Nonlinear Equations.) 1996, {\it Springer 
Verlag,}  Berlin - Heidelberg - New York. \\
11. Man Kam Kwong, Anton Zetil. Norm Inequalities for Derivatives and Differences. 1992,
{\it Lecture Notes in Mathematics.} Springer Verlag, Berlin - Heidelberg - New York. \\
12. V.V.Buldygin, D.M.Mushtary, E.I.Ostrovsky, A.I. Puchalsky. New Trends in Probability 
Theory and Statistics. 1992, {\it Springer Verlag}, Berlin - Heidelberg - New York - Amsterdam.\\
13. M.E.Taylor. Pseudodifferential Operators. 1981, {\it Princeton University Press}, Princeton,
New Jersey.\\
14. E.I.Ostrovsky. Nonlinear ODE in Orlicz spaces. Proceedings of the fourth International 
Hahn's Conference, {\it University of Chernivtzi,} 2004, p. 174 - 175.\\

%\newpage

%\begin{center}

%{\bf  APPLICATION FORM FOR PARTICIPANTS }\\
%\end{center}

%{\bf First name:} Eugene \\

%{\bf Last (family) name:} Ostrovsky \\

%{\bf Institution:} University Ben Gurion, Beer - Sheva city, Israel, Ben Gurion street, 2;
%building 58, Department of Mathematic, PO BOX 61. \\

%{\bf Location:} Israel, 56209, Rehovot city, Shkolnik street, H.5 Ap. 8; \\
%Beer - Sheva, Ben - Gurion University, Dept. of Mathematic and Computer Science.\\

%{\bf Phone numbers.} Home: (972) - (08) - 9 - 45 - 16 - 13, \  Work: (972) -
%(08) - 647 - 78 - 49. \\

%{\bf E - mail address:} galaostr@cs.bgu.ac.il; \ galo@list.ru \\

%{\bf Title of report:} Correct Solvability of Non - Linear Ordinary Differential Equations 
% in Orlicz Spaces.\\

\end{document}